\documentclass[reqno,intlimits,twoside]{article}

\usepackage{fancyhdr}
\usepackage{amsmath}
\usepackage{amssymb}
\usepackage{cite}
\usepackage{amsthm}
\usepackage{color}

\newtheorem{lemma}{Lemma}[section]
\newtheorem{theorem}{Theorem}[section]

\theoremstyle{definition}
\newtheorem{definition}{Definition}[section]

\newtheorem{remark}{Remark}[section]

\numberwithin{equation}{section}

\begin{document}

\thispagestyle{empty}

\begin{center}
    \mbox{\hskip-0.3cm \fontsize{12}{14pt}\selectfont
    \textbf{Global existence for damped $\sigma$-evolution equations with nonlocal nonlinearity}} \\[3mm]
\end{center}

\medskip

\ \!\!\!\hrulefill \
\begin{center}
\textbf{\fontsize{11}{14pt}\selectfont Khaldi Said}
\end{center}
\vskip+1cm
\thispagestyle{empty}

\noindent \textbf{Abstract.}
In this research, we would like to study the global (in time) existence of small data solutions to the following damped $\sigma$-evolution equations with nonlocal (in space) nonlinearity: 
\begin{equation*}
\partial_{t}^{2}u+(-\Delta)^{\sigma}u+\partial_{t}u+(-\Delta)^{\sigma}\partial_{t}u=I_{\alpha}(|u|^{p}), \ \  t>0, \ \ x\in \mathbb{R}^{n},
\end{equation*}
where $\sigma\geq1$, $p>1$ and  $I_{\alpha}$ is the Riesz potential of power nonlinearity $|u|^{p}$ for any $\alpha\in (0,n)$. More precisely, by using the $(L^{m}\cap L^{2})-L^{2}$ and $L^{2}-L^{2}$ linear estimates, where $m\in[1,2]$, we show the new influence of the parameter $\alpha$ on the admissible ranges of the exponent $p$.

\vskip+0.5cm

\noindent \textbf{2010 Mathematics Subject Classification.}
    42B10, 35G10, 35B45

\medskip

\noindent \textbf{Key words and phrases.}
    $\sigma$-evolution equation; global existence; Riesz potential.

\vskip+0.4cm


 \newpage
\pagestyle{fancy}
\fancyhead{}
\fancyhead[EC]{\hfill \textsf{\footnotesize Khaldi Said}}
\fancyhead[EL,OR]{\thepage}
\fancyhead[OC]{\textsf{\footnotesize Global existence for damped $\sigma$-evolution equations with nonlocal nonlinearity}\hfill }
\fancyfoot[C]{}
\renewcommand\headrulewidth{0.5pt}

\section{Introduction}

As stated in the abstract, the semi-linear damped $\sigma$-evolution models we have here in mind are: 
\begin{equation}\label{1.1}
\partial_{t}^{2}u+(-\Delta)^{\sigma}u+\partial_{t}u+(-\Delta)^{\sigma}\partial_{t}u=I_{\alpha}(\left|u\right| ^{p}), \ \ u(0,x)=u_{0}(x),\ \partial_{t}u(0,x)=u_{1}(x),
\end{equation}
with $\ \sigma \in [1,\infty), \ t\in[0,\infty)$, $x\in\mathbb{R}^{n}, \ n\geq 1$ and $I_{\alpha}$ is the Riesz potential defined by convolution product: 
\begin{equation}
(I_{\alpha}\left|u\right|^{p})(t,x)=\int_{\mathbb{R}^{n}}^{}|y-x|^{-(n-\alpha)}|u(y,t)|^{p}dy, \ \ 0<\alpha<n, \ \ p>1.
\end{equation}
\\
In (\ref{1.1}), the operator $(-\Delta)^{\sigma}$ denotes the fractional Laplacian with symbol $|\xi|^{2\sigma}$, i.e., 
$$
\mathcal{F}\left( (-\Delta)^{\sigma}f\right) = |\xi|^{2\sigma}\mathcal{F}\left(f\right) (\xi), \ \  \xi\in\mathbb{R}^{n}, \ \  |\xi|=(\xi_{1}^{2}+\cdots+\xi_{n}^{2})^{1/2},
$$ 
where $\mathcal{F}$ is the Fourier transform. In this research, for the sake of simplicity, we will choose
\begin{equation}\label{1.2}
u(0,x)=0 \ \ \textit{and} \ \  u_{1}\in L^{2}(\mathbb{R}^{n})\cap L^{m}(\mathbb{R}^{n}), \ m\in [1,2).
\end{equation}

By the recent paper \cite{DaoMichihisa}, the authors found the following critical exponent $p_{crit}(n,m,\sigma)$ of the Cauchy problem (\ref{1.1}) with the usual power nonlinearity $|u|^{p}$ in the right hand side 
\begin{equation}\label{1.3}
p_{crit}(n,m,\sigma)=1+\frac{2m\sigma}{n}, \ \ m\in[1,2).
\end{equation}

This exponent guarantees both the global (in time) existence of small data Sobolev solutions for $p > p_{crit}(n,m,\sigma)$, and blow-up in finite time for $1 < p <p_{crit}(n,m,\sigma)$. 

Fortunately, there are two ways to catch the critical exponents of some semi-linear evolution equations. The first one is by using some \textit{sharp} linear estimates and Banach fixed point theorem in the prove of global existence result, while, the second one is to prove the blow-up result by the test function method based on contradiction argument.

Thanks to the presence of two damping mechanism in (\ref{1.1}), the solutions of the linear equation 
\begin{equation}\label{2.1}
\partial_{t}^{2}u+(-\Delta)^{\sigma}u+\partial_{t}u+(-\Delta)^{\sigma}\partial_{t}u=0, \ u(0,x)=0,\   \partial_{t}u(0,x)=u_{1}(x), 
\end{equation}
satisfies some sharp $(L^{m}\cap L^{2})-L^{2}$ and $L^{2}-L^{2}$ linear estimates (Lemma \ref{LE1}). More precisely, the $L^{2}-L^{2}$ estimates can be seen as the $(L^{m}\cap L^{2})-L^{2}$ estimates after formally set $m=2$. This fact has some benefits in the treatment of Duhamel's integral.

So, in this research we will rely on the first way, which is the use of these sharp linear estimates together with Banach fixed point theorem and prove not only the global (in time) existence of small data solutions to (\ref{1.1}), but also study the influence of the parameter $\alpha \in (0,n)$ on the critical exponent (\ref{1.3}) and some admissible ranges of exponent $p$. 

The paper is structured as follows: In Section \ref{Main tools} we have some preliminaries which play an essential role to prove our main result in Section \ref{Main results}. 
\section{Preliminaries }\label{Main tools}
First, we have the following notations.
\begin{itemize}
	\item we write $f\lesssim g$ when there exists a constant $c>0$ such that $f \leq cg$;
	\item $L^{m}(\mathbb{R}^{n})$ denote the usual Lebesgue spaces, with $m\in [1,2)$; 
	\item $H^{\sigma}(\mathbb{R}^{n})$ mean Sobolev spaces as defined below (see \cite{EbertReissig} p 445):
	$$H^{\sigma}(\mathbb{R}^{n}):=\left\lbrace f\in S'(\mathbb{R}^{n}): \|f\|_{H^{\sigma}(\mathbb{R}^{n})}=\|(1+|\cdot|^{2})^{\frac{\sigma}{2}}\mathcal{F}(f)\|_{L^{2}(\mathbb{R}^{n})}<\infty\right\rbrace.$$
\end{itemize}

We now introduce the Riesz potential and the boundedness property from $L^{q}$ to $L^{r}$ spaces.
\begin{definition}\cite[p117]{stien}
	Let $\alpha \in(0,n)$. We formally define the normalized Riesz potential as:
	$$(I_{\alpha}f)(x):=\frac{\Gamma((n-\alpha)/2)}{\pi^{n/2}2^{\alpha}\Gamma(\alpha/2)} \int_{\mathbb{R}^{n}}^{}\frac{f(y)}{|y-x|^{n-\alpha}}dy, \ \ x\in \mathbb{R}^{n}.$$
\end{definition}
\begin{lemma}\cite[p119]{stien}
	If $f \in L^{q}$ for some $q\in(1,n/\alpha)$, then  $I_{\alpha}f \in L^{r}$, where
	$$\|I_{\alpha}f\|_{L^{r}}\lesssim \|f\|_{L^{q}}, \ \ \frac{1}{q}-\frac{1}{r}=\frac{\alpha}{n}.$$	
\end{lemma}
We now recall the following sharp linear estimates which is very important tool to demonstrate Theorem \ref{GlobalExistence1}.
\begin{lemma}\cite[Proposition 2.1]{DaoMichihisa}\label{LE1}
	Let $m\in [1,2)$. Then, the Sobolev solutions $u$ to the linear equation (\ref{2.1}) satisfy the $(L^{m}\cap L^{2})-L^{2}$ estimates:
	\begin{align} 
	\|\partial_{t}^{j}(-\Delta)^{a/2}u(t,\cdot)\|_{L^{2}} &\lesssim (1+t)^{-\frac{n}{2\sigma}\left(\frac{1}{m}-\frac{1}{2}\right)-\frac{a}{2\sigma}-j}\|u_{1}\|_{(L^{m}\cap H^{[a+2(j-1)\sigma]^{+}})}, \label{2.2} 
	\end{align}
	and the $L^{2}-L^{2}$ estimates:
	\begin{equation}
	\|\partial_{t}^{j}(-\Delta)^{a/2}u(t,\cdot)\|_{L^{2}} \lesssim (1+t)^{-\frac{a}{2\sigma}-j}\|u_{1}\|_{H^{[a+2(j-1)\sigma]^{+}}}, \label{2.3}
	\end{equation}
	for any $a\geq 0$, $j=0,1$ and for all space dimensions $n\geq 1$, where $[\cdot]^{+}=\max\{0,\cdot\}$.
\end{lemma}
\smallskip
We also need to the fractional Gagliardo-Nirenberg inequality in the following lemma.
\begin{lemma} [\cite{HajaiejMolinetOzawa, EbertReissig}]\label{FGN}
	Let $1<q<\infty$, $\sigma>0$. Then, the following fractional Gagliardo-Nirenberg inequality holds for all $y\in H^{\sigma}(\mathbb{R}^{n})$
	\[\|y\|_{L^{q}(\mathbb{R}^{n})}\lesssim \|(-\Delta)^{\sigma/2}y\|_{L^{2}(\mathbb{R}^{n})}^{\theta_{q}}\,\|y\|_{L^{2}(\mathbb{R}^{n})}^{1-\theta_{q}},\]
	where \[\theta_{q}=\frac{n}{\sigma}\left(\frac{1}{2}-\frac{1}{q}\right)\in\left[0,1\right].\]
\end{lemma}
Finally, this integral inequality is used to deal with the Duhamel's integral. 
\begin{lemma}[\cite{EbertReissig}]\label{Integral inequality}
	Let $a, b\in\mathbb{R}$ such that $\max\{a,b\}>1$. Then, it holds
	$$\int_{0}^{t}(1+t-\tau)^{-a}(1+\tau)^{-b}d\tau\lesssim (1+t)^{-\min\{a,b\}}.$$ 
\end{lemma} 

\smallskip
\section{Main Result}\label{Main results}
Our main result is reads as follows.
\begin{theorem}\label{GlobalExistence1}
	Let us consider the Cauchy problem (\ref{1.1}) with $\sigma\geq1$, $p>1$ and $\alpha\in (0,n)$. Let $m \in [1,2)$, we assume the following conditions for $p$ and the dimension $n$:
	\begin{equation}\label{1.4}
	\left\lbrace  
	\begin{matrix}
	\frac{2}{m}+\frac{2\alpha}{n}\leq p \leq \frac{n+2\alpha}{n-2\sigma} \hfill& {if }\ \ 2\sigma<n\leq \frac{4\sigma+\sqrt{16\sigma(\sigma+m(2-m)\alpha)}}{2(2-m)}, 
	&\cr
	\\
	\frac{2}{m}+\frac{2\alpha}{n}\leq p < \infty \hfill & {if } \ \
	1\leq n \leq 2\sigma .
	&\cr
	\end{matrix}\right.  
	\end{equation}
	Moreover, we suppose
	\begin{equation}\label{1.5}  
	p>1+\frac{(2\sigma+\alpha)m}{n}.
	\end{equation}\\
	Then, there exists a constant $\varepsilon_0>0$ such that for any data 
	$$u_{1}\in L^{2}(\mathbb{R}^{n})\cap L^{m}(\mathbb{R}^{n}) \ \textit{with}\  \left\|u_{1}\right\|_{L^{2}(\mathbb{R}^{n})\cap L^{m}(\mathbb{R}^{n})}<\varepsilon_0,$$ 
	we have a uniquely determined globally (in time) solution
	$$u\in\mathcal{C}\left([0,\infty), H^{\sigma}(\mathbb{R}^{n})\right)\cap \mathcal{C}^{1}\left([0,\infty),L^{2}(\mathbb{R}^{n})\right)$$
	to (\ref{1.1}). Furthermore, the solution satisfies the estimates:
	\begin{equation*}
	\left\|u(t,\cdot)\right\|_{L^{2}(\mathbb{R}^{n})} \lesssim (1+t)^{-\frac{n}{2\sigma}\left(\frac{1}{m}-\frac{1}{2}\right)}\left\|u_{1}\right\|_{L^{2}(\mathbb{R}^{n})\cap L^{m}(\mathbb{R}^{n})},
	\end{equation*}
	\begin{equation*}
	\left\|\partial_{t}u(t,\cdot)\right\|_{L^{2}(\mathbb{R}^{n})} \lesssim (1+t)^{-\frac{n}{2\sigma}\left(\frac{1}{m}-\frac{1}{2}\right)-1}\left\|u_{1}\right\|_{L^{2}(\mathbb{R}^{n})\cap L^{m}(\mathbb{R}^{n})},
	\end{equation*}
	\begin{equation*} 
	\left\|(-\Delta)^{\sigma/2}u(t,\cdot)\right\|_{L^{2}(\mathbb{R}^{n})} \lesssim (1+t)^{-\frac{n}{2\sigma}\left(\frac{1}{m}-\frac{1}{2}\right)-\frac{1}{2}}\left\|u_{1}\right\|_{L^{2}(\mathbb{R}^{n})\cap L^{m}(\mathbb{R}^{n})}.
	\end{equation*}
\end{theorem}
\begin{remark}
	As it happens with the power nonlinearity $|u|^{p}$, the conditions (\ref{1.5}) is assumed to get the same decay estimates of the semi-linear model with those of the corresponding linear model (\ref{2.1}). The bounds (\ref{1.4}) on $p$ and $n$ appear due to the application of Gagliardo-Nirenberg inequality from Lemma (\ref{FGN}) and the boundedness property of Riesz potential.
\end{remark}
\begin{remark} Theorem \ref{GlobalExistence1} showed the influence of the Riesz potential not only on the critical exponent but also on the admissible ranges of exponent $p$. It is clear that if we let $\alpha \to 0$ then our results in Theorems \ref{GlobalExistence1} tend to be exactly the same as in the cited paper \cite{DaoMichihisa}. This potential has no influences on the decay estimates, because it is nonlocal in space which is in contrast to the effect of nonlinear memory (nonlocal in time)
	$$\int_{0}^{t}(t-s)^{-\gamma}|u(s,x)|^{p}ds, \ \ \gamma \in (0,1),  \ \ p>1.$$  
\end{remark}

Since we are dealing with semi-linear Cauchy problems, we use the Banach's fixed point theorem inspired from the book \cite[Page 303]{EbertReissig}. Here, we need to define a family of evolution spaces $X(T)$ for any $T>0$ with suitable norm $\|\cdot\|_{X(T)}$, also, we define an operator $$O :u\in X(T)\longmapsto Ou(t,x)=G_{\sigma}(t,x)\ast u_{1}(x)+  \int_{0}^{t}G_{\sigma}(t-\tau,x)\ast I_{\alpha}(|u(\tau,x)|^{p})d\tau.$$
If this operator satisfies the two inequalities: 
\begin{align}
\|Ou\|_{X(T)} &\lesssim \left\|u_{1}\right\|_{L^{2}(\mathbb{R}^{n})\cap L^{m}(\mathbb{R}^{n})}+\|u\|_{X(T)}^{p}, \ \  \forall u \in X(T), \label{3.3} \\ 
\|Ou-Ov\|_{X(T)} &\lesssim \|u-v\|_{B(T)}\Big( \|u\|_{X(T)}^{p-1}+\|v\|_{X(T)}^{p-1}\Big), \ \  \forall u, v \in X(T), \label{3.4} 
\end{align} 
then, one can deduce the existence and uniqueness of a global (in time) solutions of (\ref{1.1}). Here, the smallness of the initial data $\left\|u_{1}\right\|_{L^{2}(\mathbb{R}^{n})\cap L^{m}(\mathbb{R}^{n})}<\varepsilon_0$ imply that the operator $O$ maps balls of $X(T)$ into balls of $X(T)$.

The Banach space $X(T)$ is defined for all $T>0$ as follows: 
$$
X(T):=\mathcal{C}\left([0,T], H^{\sigma}\right)\cap \mathcal{C}^{1}\left([0,T],L^{2}\right).
$$
We equip this space with the norm 
\begin{align}
\|u\|_{X(T)}&=\sup_{0\leq t\leq T}\Big((1+t)^{\frac{n}{2\sigma}\left(\frac{1}{m}-\frac{1}{2}\right)}\|u(t,\cdot)\|_{L^{2}}+(1+t)^{ \frac{n}{2\sigma}\left(\frac{1}{m}-\frac{1}{2}\right)+\frac{1}{2}}\|(-\Delta)^{\sigma/2}u(t,\cdot)\|_{L^{2}} \nonumber \\
&\hspace{3cm} 
+(1+t)^{ \frac{n}{2\sigma}\left(\frac{1}{m}-\frac{1}{2}\right)+1}\|\partial_{t}u(t,\cdot)\|_{L^{2}}\Big). \label{3.2}
\end{align} 

\begin{proof}
	\textit{Step 1:} By using the linear estimates, it is clear that the function $$u^{L}(t,x)=G_{\sigma}(t,x)\ast u_{1}(x)$$ is belongs to $X(T)$ and we have:
	\begin{align*}
	\|u^{L}\|_{X(T)}&=\sup_{0\leq t\leq T}\Big((1+t)^{\frac{n}{2\sigma}\left(\frac{1}{m}-\frac{1}{2}\right)}\|u^{L}(t,\cdot)\|_{L^{2}}+(1+t)^{ \frac{n}{2\sigma}\left(\frac{1}{m}-\frac{1}{2}\right)+\frac{1}{2}}\|(-\Delta)^{\sigma/2}u^{L}(t,\cdot)\|_{L^{2}} \nonumber \\
	&\hspace{3cm} 
	+(1+t)^{ \frac{n}{2\sigma}\left(\frac{1}{m}-\frac{1}{2}\right)+1}\|\partial_{t}u^{L}(t,\cdot)\|_{L^{2}}\Big)\lesssim \left\|u_{1}\right\|_{L^{2}(\mathbb{R}^{n})\cap L^{m}(\mathbb{R}^{n})}.
	\end{align*} 
	\textit{Step 2:} To conclude inequality (\ref{3.3}), we prove
	\begin{equation}\label{3.6}
	\|u^{N}\|_{B(T)}\lesssim \|u\|_{X(T)}^{p}.
	\end{equation}
	We divide the interval $[0,t]$ into two sub-intervals $[0,t/2]$ and $[t/2,t]$ where we use the $L^{m}-L^{2}$ linear estimates if $\tau\in[0,t/2]$ and $L^{2}-L^{2}$ estimates if $\tau\in[t/2,t]$. From Lemma \ref{LE1} we estimate:
	\begin{align}
	\|u^{N}(t,\cdot)\|_{L^{2}}&\lesssim\int_{0}^{t/2}(1+t-\tau)^{-\frac{n}{2\sigma}\left(\frac{1}{m}-\frac{1}{2}\right)}\left( \left\|I_{\alpha}(|u|^{p})(\tau,\cdot)\right\| _{L^{2}}+\left\|I_{\alpha}(|u|^{p})(\tau,\cdot)\right\| _{L^{m}}\right) d\tau \nonumber \\
	&\hspace{3cm} 
	+\int_{t/2}^{t}\left\|I_{\alpha}(|u|^{p})(\tau,\cdot)\right\|_{L^{2}}d\tau,  \label{3.7}
	\end{align}
	\begin{align}
	\|\partial_{t}u^{N}(t,\cdot)\|_{L^{2}}&\lesssim\int_{0}^{t/2}(1+t-\tau)^{-\frac{n}{2\sigma}\left(\frac{1}{m}-\frac{1}{2}\right)-1}\left( \left\|I_{\alpha}(|u|^{p})(\tau,\cdot)\right\| _{L^{2}}+\left\|I_{\alpha}(|u|^{p})(\tau,\cdot)\right\|_{L^{m}}\right)d\tau \nonumber \\
	&\hspace{3cm} 
	+\int_{t/2}^{t}(1+t-\tau)^{-1}\left\|I_{\alpha}(|u|^{p})(\tau,\cdot)\right\|_{L^{2}}d\tau,  \label{3.8}
	\end{align}
	\begin{align}
	\|(-\Delta)^{\sigma/2}u^{N}(t,\cdot)\|_{L^{2}}&\lesssim\int_{0}^{t/2}(1+t-\tau)^{-\frac{n}{2\sigma}\left(\frac{1}{m}-\frac{1}{2}\right)-\frac{1}{2}}\left( \left\|I_{\alpha}(|u|^{p})(\tau,\cdot)\right\| _{L^{2}}+\left\|I_{\alpha}(|u|^{p})(\tau,\cdot)\right\|_{L^{m}}\right)d\tau \nonumber \\
	&\hspace{3cm} 
	+\int_{t/2}^{t}(1+t-\tau)^{-\frac{1}{2}}\left\|I_{\alpha}(|u|^{p})(\tau,\cdot)\right\|_{L^{2}}d\tau.  \label{3.9}
	\end{align}
	Now, to control these integrals, we use the boundedness property of Riesz potential as follows:
	$$\|I_{\alpha}(|u|^{p})(\tau,\cdot)\|_{L^{2}}\lesssim \|u(\tau,\cdot)\|^{p}_{L^{\frac{2np}{n+2\alpha}}}, \ \ \|I_{\alpha}(|u|^{p})(\tau,\cdot)\|_{L^{m}}\lesssim \|u(\tau,\cdot)\|^{p}_{L^{\frac{mnp}{n+m\alpha}}}.$$ 
	By the fractional Gagliardo-Nirenberg inequality from Lemma \ref{FGN}, we can estimate these later norms: 
	$$\|u(\tau,\cdot)\|^{p}_{L^{\frac{2np}{n+2\alpha}}},\ \ \|u(\tau,\cdot)\|^{p}_{L^{\frac{mnp}{n+m\alpha}}}.$$
	Here, we know from (\ref{3.2}):
	\begin{equation*}
	(1+\tau)^{\frac{n}{2\sigma}\left(\frac{1}{m}-\frac{1}{2}\right)+\frac{1}{2}}\|(-\Delta)^{\sigma/2}u(\tau,\cdot)\|_{L^{2}}\lesssim\|u\|_{X(T)},
	\end{equation*}
	\begin{equation*}
	(1+\tau)^{\frac{n}{2\sigma}\left(\frac{1}{m}-\frac{1}{2}\right)}\|u(\tau,\cdot)\|_{L^{2}}\lesssim\|u\|_{X(T)}.
	\end{equation*}
	So, we can estimate the above norms as follows:
	\begin{equation}
	\left\|u(\tau,\cdot)\right\| _{L^{\frac{snp}{n+s\alpha}}}^{p}\lesssim(1+\tau)^{-\frac{np}{2m\sigma}+\frac{n}{2\sigma}\left(\frac{1}{s}+\frac{\alpha}{n}\right) }\|u\|_{X(T)}^{p}, \ s=m, 2, \label{3.10}
	\end{equation}
	provided that the conditions (\ref{1.4}) are satisfied for $p$ and $n$. Hence, we conclude 
	\begin{equation}\label{3.11}
	\left\|u(\tau,\cdot)\right\|^{p}_{L^{\frac{2np}{n+2\alpha}}}+\left\|u(\tau,\cdot)\right\|^{p}_{L^{\frac{mnp}{n+m\alpha}}}\lesssim (1+\tau)^{-\frac{np}{2m\sigma}+\frac{n}{2\sigma}\left(\frac{1}{m}+\frac{\alpha}{n}\right) }\|u\|_{X(T)}^{p}. 
	\end{equation}
	By these equivalences:
	$$(1+t-\tau)\approx (1+t) \ \text{if} \ \tau \in[0,t/2],\  \ (1+\tau)\approx (1+t) \ \text{if} \ \tau \in[t/2,t] $$
	and Lemma \ref{Integral inequality}, we estimates the first integral of $u^{N}$ over $[0,t/2]$ as follows:
	$$
	\int_{0}^{t/2}(1+t-\tau)^{-\frac{n}{2\sigma}\left(\frac{1}{m}-\frac{1}{2}\right)}(1+\tau)^{-\frac{np}{2m\sigma}+\frac{n}{2\sigma}\left(\frac{1}{m}+\frac{\alpha}{n}\right) }\|u\|_{X(T)}^{p}d\tau$$
	$$
	\lesssim(1+t)^{-\frac{n}{2\sigma}\left(\frac{1}{m}-\frac{1}{2}\right)}\|u\|_{X(T)}^{p}\int_{0}^{t/2}(1+\tau)^{-\frac{np}{2m\sigma}+\frac{n}{2\sigma}\left(\frac{1}{m}+\frac{\alpha}{n}\right) }d\tau
	$$
	$$
	\lesssim(1+t)^{-\frac{n}{2\sigma}\left(\frac{1}{m}-\frac{1}{2}\right)}\|u\|_{X(T)}^{p}
	$$
	provided that $p>1+\frac{(2\sigma+\alpha)m}{n}$. For the second integral over $[t/2,t]$ we have:
	
	$$\int_{t/2}^{t}(1+\tau)^{-\frac{np}{2m\sigma}+\frac{n}{2\sigma}\left(\frac{1}{2}+\frac{\alpha}{n}\right) }\|u\|_{X(T)}^{p}d\tau\lesssim(1+t)^{1-\frac{np}{2m\sigma}+\frac{n}{2\sigma}\left(\frac{1}{2}+\frac{\alpha}{n}\right) }\|u\|_{X(T)}^{p}.$$
	Thanks to $p>1+\frac{(2\sigma+\alpha)m}{n}$, we arrive to the desired estimate for $u^{N}$ 
	$$
	(1+t)^{\frac{n}{2\sigma}\left(\frac{1}{m}-\frac{1}{2}\right)}\|u^{N}(t,\cdot)\|_{L^{2}}\lesssim\|u\|_{X(T)}^{p}.
	$$
	Without repetition, we can prove again:
	$$ 
	(1+t)^{\frac{n}{2\sigma}\left(\frac{1}{m}-\frac{1}{2}\right)-1}\|\partial_{t}u^{N}(t,\cdot)\|_{L^{2}}\lesssim\|u\|_{X(T)}^{p},
	$$
	$$
	(1+t)^{\frac{n}{2\sigma}\left(\frac{1}{m}-\frac{1}{2}\right)-\frac{1}{2}}\|(-\Delta)^{\sigma/2}u^{N}(t,\cdot)\|_{L^{2}}\lesssim\|u\|_{X(T)}^{p},
	$$ 
	and the inequality (\ref{3.6}) is now proved, that is (\ref{3.3}).\\
	\textit{Step 3:} To prove (\ref{3.4}) we choose two elements $u$, $v$ belong to $X(T)$, and we write
	$$Ou-Ov=\int_{0}^{t}G_{\sigma}(t-\tau,x)\ast I_{\alpha}(|u(\tau,x)|^{p}-|v(\tau,x)|^{p})d\tau.$$
	We divide all the integrals as above, we write again:
	\begin{eqnarray}
	\nonumber 	\|(u^{N}-v^{N})(t,\cdot)\|_{L^{2}}&\lesssim&\int_{0}^{t/2}(1+t-\tau)^{-\frac{n}{2\sigma}\left(\frac{1}{m_{2}}-\frac{1}{2}\right)}\left\|I_{\alpha}(|u(\tau,x)|^{p}-|v(\tau,x)|^{p})\right\|_{L^{m}\cap L^{2}}d\tau\\
	&&\qquad\qquad\qquad	+\int_{t/2}^{t}\left\|I_{\alpha}(|u(\tau,x)|^{p}-|v(\tau,x)|^{p})\right\| _{L^{2}}d\tau \nonumber,\\
	\nonumber	\|\partial_{t}(u^{N}-v^{N})(t,\cdot)\|_{L^{2}}&\lesssim&\int_{0}^{t/2}(1+t-\tau)^{-\frac{n}{2\sigma}\left(\frac{1}{m}-\frac{1}{2}\right)-1}\left\|I_{\alpha}(|u(\tau,x)|^{p}-|v(\tau,x)|^{p})\right\| _{L^{m}\cap L^{2}}d\tau\\
	&&\qquad+\int_{t/2}^{t}(1+t-\tau)^{-1}\left\|I_{\alpha}(|u(\tau,x)|^{p}-|v(\tau,x)|^{p})\right\|_{L^{2}}d\tau,\nonumber \\
	\nonumber	\|(-\Delta)^{\sigma/2}(u^{N}-v^{N})(t,\cdot)\|_{L^{2}}&\lesssim&\int_{0}^{t/2}(1+t-\tau)^{-\frac{n}{2\sigma}\left(\frac{1}{m}-\frac{1}{2}\right)-\frac{1}{2}}\left\|I_{\alpha}(|u(\tau,x)|^{p}-|v(\tau,x)|^{p})\right\| _{L^{m}\cap L^{2}}d\tau\\
	&&
	\qquad+\int_{t/2}^{t}(1+t-\tau)^{-\frac{1}{2}}\left\|I_{\alpha}(|u(\tau,x)|^{p}-|v(\tau,x)|^{p})\right\|_{L^{2}}d\tau.\nonumber
	\end{eqnarray}
	Thanks to the boundedness property of Reisz potential, we have 
	$$\|I_{\alpha}(|u(\tau,x)|^{p}-|v(\tau,x)|^{p})\|_{L^{2}}\lesssim \||u(\tau,x)|^{p}-|v(\tau,x)|^{p}\|_{L^{\frac{2n}{n+2\alpha}}},$$
	$$ \|I_{\alpha}(|u(\tau,x)|^{p}-|v(\tau,x)|^{p})\|_{L^{m}}\lesssim \||u(\tau,x)|^{p}-|v(\tau,x)|^{p}\|_{L^{\frac{mn}{n+m\alpha}}}.$$ 
	Now, by the H\"{o}lder's inequality, we derive  for $r=\frac{mn}{n+m\alpha}, \frac{2n}{n+2\alpha}$ the following
	\begin{equation}\label{5.23}
	\| |u(\tau,\cdot)|^{p}-|v(\tau,\cdot)|^{p}\|_{L^{r}}\leq \|u(\tau,\cdot)-v(\tau,\cdot)\|_{L^{r p}}\left(\|u(\tau,\cdot)\|_{L^{r p}}^{p-1}+\|v(\tau,\cdot)\|_{L^{rp}}^{p-1} \right).
	\end{equation}
	Using the definition of the norm $\|u-v\|_{X(T)}$ and fractional Gagliardo-Nirenberg inequality we can prove the inequality (\ref{3.4}). Hence, Theorem \ref{GlobalExistence1} is proved.	
\end{proof}

\vskip+0.5cm

\noindent \textbf{Author address:}

\vskip+0.3cm

\noindent \textbf{Khaldi Said}

Laboratory of Analysis and Control of PDEs, Djillali Liabes University,  P.O. Box 89, Sidi-Bel-Abb\`{e}s 22000, Algeria

{\itshape E-mails:} \texttt{saidookhaldi@gmail.com, said.khaldi@univ-sba.dz}

\end{document}